%% file: geometric_regulation.tex
\documentclass[12pt]{article}

 \pdfoutput=1
 
\usepackage{graphicx}
\usepackage{amssymb}

  \usepackage{amsmath,amsthm,amssymb,euscript,cite}
\usepackage{graphicx,paralist,geometry}
\input dsgheader20

\allowdisplaybreaks

\newtheorem{theorem}{Theorem}[section]
\newtheorem{proposition}{Proposition}[section]
\newtheorem{definition}{Definition}[section]

\newtheorem{assumption}{Assumption}[section]
\newtheorem{remark}{Remark}[section]
\newtheorem{problem}{Problem}[section]

\numberwithin{equation}{section}

\usepackage{xcolor}

\begin{document}
 
  \begin{center}
 {\bf \Large  Approximation Methods for Geometric Regulation }\\[2ex]
 {\bf \large E. Aulisa,    \ \ D.S. Gilliam }
 \end{center}
 
 \begin{abstract}                

In these notes we collect some results from several of the authors' works in order to make available a single source and show how the approximate geometric methods for regulation have been developed, and how the control design strategy has evolved from the theoretical methods, involving  the regulator equations, to what we now call the regularized controller. In between these two extremes we developed, in  a series of works, a fairly rigorous analysis of the regularization scheme leading to the regularized dynamic regulator equations and an iterative scheme that produces very accurate tracking and disturbance rejection control laws. In our most recent work we have extended dynamic regulator equations to what we now refer to as the regularized controller. This new formulation has only recently being applied to examples including linear and nonlinear delay equations. 
\end{abstract}

\section{Geometric Regulation for Nonlinear Systems: Regulator Equations}\label{sec_regeq}
In these notes, in order to avoid a number of technical issues, we consider a Single Input Single Output  (SISO)   abstract  nonlinear control system
\begin{align}
\deriv{z}{t} &=  A z +f(z) +\Bin u +\Bd Pw, \label{eq1} \\
   z(0) &=z_0 ,  \label{eq3} \\ 
   y &=C z .  \label{eq4} 
\end{align}

In this section we will assume that the disturbance $d$ and the signal to be tracked $r$ are generated by a linear finite dimensional exogenous system (in later sections we will consider more general signals)
 \begin{align}
   \dot w & =  S w ,\ \  w(0)   = w_0, \label{eq2} \\
    d &= Pw , \ \ 
    r = Qw, \nonumber
    \end{align}
in  the state  space $\sW$, with $\dim (\sW)=N$. Here  $S$ is assumed to be neutrally stable (with respect to the origin, i.e. Lyapunov stable but not attractive). We have $Q\in \sL(\sW , \sY)$ and $P\in \sL(\sW , \sD)$. Here, $\sY$ and $\sD$ are the output and the disturbance space, respectively.

We define the system error to be
\begin{equation}
\label{e}
e(t) = r(t) - y(t)    = Q w(t) - C z(t).
\end{equation}

The main theoretical Asymptotic  Regulation  problem found in the literature is  described in  Problem \ref{prob1}.
 \begin{problem}[\bf  State Feedback  Regulation Problem]\label{prob1}
{\rm  \mbox{}\\
  Find a time dependent control law $u \in C_b(\bbr^+,\sU)$,  for the system \rr{eq1}-\rr{eq2},   so  that the tracking error, defined as    \begin{equation}
\label{Et}
e(t) = r(t)-y(t),
\end{equation}
satisfies 
 \begin{equation}
\label{er10}
  \lim _ {t \to \infty} \| e(t) \|_{\sY}  = 0\ \ \ \ \ \ \  \text{  (asymptotic regulation),  }  
\end{equation}
 for all (sufficiently small) initial conditions $w_0$ in \rr{eq2} and $z_0$ in \rr{eq3}.
 }
\end{problem}
The material in this section mainly derives from the work \cite{bg_acc08}.
\begin{assumption}\label{ass1}{\rm
\begin{enumerate}
\item The operator $A$ is  a sectorial operator,    with compact resolvent
  $R(\l,A) = (\l I-A)^{-1}$  for $\l\in \r(A)$ (the resolvent set of $A$)  and is assumed to generate  an exponentially stable   analytic semigroup $e^{At}$ in  $\sZ$.    
 In particular there exist numbers $\o>0$ and $M\ge1$ such that  
 \begin{equation}
\label{Aexpstab}
 \| \exp^{At}   \| \leq M \exp^{-\o t}. 
\end{equation}
Here the norm on the left is the operator norm. The symbol $\sL(W_1,W_2)$  denotes the set of all bounded linear operators from a  Hilbert space $W_1$ to a Hilbert space $W_2$. When $W=W_1=W_2$ we write $\sL(W)$.

The operator $A$ generates an infinite scale of  Banach spaces denoted by $\sZ^{s}$.  For $s>0$, let $\sZ^s =  \ran((-A)^{-s})$ with norm $ \|\vp\|_s = \|(-A)^{s}\vp\|$, where, as usual, the norm on the right, $\|\cdot\|$, denotes the norm in $\sZ$.  For $s>0$ the space $\sZ^{-s}$ is the completion of $\sZ$ with norm  $ \|\vp\|_{-s} = \|(-A)^{-s}\vp\|$.
For a detailed discussion of  fractional powers of sectorial operators and properties of the scale of  Banach spaces, see \cite{henry,kato,pazy,Engel,haase}. 
  In this work we follow the  definition  of sectorial operator used in \cite{Engel} which is different from the one found, for example, in \cite{henry}. In particular, our $A$ would be $-A$ in \cite{henry}. Assumption \ref{ass1}  implies that  the spectrum of $A$ consists of discrete eigenvalues whose only limit point is at infinity and  lies in a sector in the left half complex plane.

 \item The input space is  $\sU= \bbr$ and   $\Bin \in \sL(\sU, \sZ^{-\sb})$ for some $0\leq \sb <1$.  
For a given input vector $v(t) \in \sU$,
it follows
$$\Bin v(t)=  b v   \text{  for some  } b  \in \sZ^{-\sb}.$$
For the operator norm the following notation is used
$$   \|\Bin \|_{\sL(\sU, \szmb)} \defeq \|\Bin \|_{-s_b} =      \|\sa{-s_b} b \|   .$$
In the work \cite{bg_acc08} it is assumed that $\Bin$ is bounded in the Hilbert state space so that $\sb=0$. Therefore, for the results in this section we assume $b\in \sZ$.

\item The disturbance space is  $\sD= \bbr$  and  $\Bd \in \sL(\sD, \sZ^{-\sd})$  for some $0\leq \sd <1$.  For given disturbance  
$d(t)  \in  \sD  $,  it follows
$$\Bd d(t)=  b_d d (t) \text{  for some } {b_d}  \in \sZ^{-\sd}.$$
For the operator norm the following notation  is used
$$   \|\Bd\|_{\sL(\sD, \szmd)} \defeq \|\Bd\|_{-s_d} =       \|\sa{-s_d} b_{d} \|   .$$
 Once again, in this section we assume $\Bd$ is bounded so that $\sd=0$ and $b_d\in \sZ$.
 
\item The output space is  $\sY= \bbr $ and   $C \in \sL(\sZ^{\sc},\sY)$ for some $0\leq \sc <1$. Therefore, for $\vp\in \sZ^{s_c}$  
$$\|C\vp\|_{\sY} = \|C\sa{-s_c } \sa{s_c } \vp\|_{\sY} \leq \|C\sa{-s_c }\|_{\sY} \, \|\sa{s_c } \vp\|_{\sZ} =\|C\sa{-s_c }\|_{\sL(\sZ,\sY)}\|\vp\|_{s_c},$$
which implies
$$ \| C\|_{ \sL(\sZ^{s_c}, \sY) } \leq  \|C\sa{-s_c }\|_{\sL(\sZ,\sY)}.$$
Again the following notation is used 
$$ \| C\|_{ \sL(\sZ^{s_c}, \sY) }  \defeq \|C\|_{-s_c}.$$ 
 
 \item $f\, :\,  \sZ \ra \sZ $ is a smooth nonlinear function so that, in particular,  for all $\vp_1$, $\vp_2 \in  \sZ$
  \begin{equation}
\label{eq_flip}
\| f(\vp_1) - f(\vp_2)\|_{\sZ} \leq \e \|\vp_1-\vp_2\|_{\a}.
\end{equation}
for some $0<\a<1$ and $\e>0$.
\item We assume that the exosystem has the origin as a  neutrally stable equilibrium, i.e.,   $w=0$ is a fixed point which  is Lyapunov stable but not attracting. A center is an example of such a fixed point. This, in particular, implies   $\s(S)\subset i\bbr$ (i.e., the spectrum of $S$  is on the imaginary axis) and  has no non-trivial Jordan blocks.  

\end{enumerate}
}
\end{assumption}

 \begin{remark}\label{regular}{\rm
 Referring to the notation and terminology  in \cite{staffans} and using Assumptions \ref{ass1}   with $s_c+s_b<1$ and $s_c+s_d<1$, Theorem  5.7.3 in \cite{staffans} tells us that the linearization of system \rr{eq1}--\rr{eq2} is an $L^1$-well-posed  linear system in ($\sY,\sZ^{\g},\sU)$ for $\g\in(s_c-1,-s_b]$  and it is also uniformly regular.  A detailed discussion of these points is not needed in this paper, however it precisely corresponds to the $L^1$-well-posedness that motivates the bounds obtained in the
 error analysis of the proposed iterative scheme discussed below.
 }
 \end{remark}

With the above, the main result from \cite{bg_acc08} is 
\smallskip
 \begin{theorem} \label{thm1}     Under  Assumption \ref{ass1},  the  (local)  State Feedback Regulator Problem \ref{prob1} for \rr{eq1}-\rr{eq4} 
is solvable if, and only if,    there exist mappings $ \pi \, :\, \sW\ra D(A)\subset \sZ$   and  $ \g\, :\, \sW\ra \sY$    satisfying the ``regulator equations,''  
\begin{align}
\label{sfrp1}
&  \pderiv{\pi}{w} Sw  = A\pi(w) + f(\pi(w))+\Bin\g(w) + \Bd Pw, \\
\label{sfrp2}
& C \pi(w) =Qw .
\end{align}
  In this case a feedback law solving the state feedback regulator problem is given by 
\begin{equation} u(t) = \g(w)(t).\label{u} \end{equation}

\end{theorem}
\smallskip
Let us introduce the linearized version of $\g$
$$u  = \g(w) = \G w  +\wt{\g}(w),$$
$$ \G\in \call(\sW,\sU), \ \ \wt{\g}(0)=0,\   \pderiv{\wh{\g}}{w}(0) =0.$$

Modulo the inherent  technical difficulties that arise in infinite dimensions, Theorem \ref{thm1} can be obtained using an argument similar to that given in \cite{IB}. Indeed,  under the assumptions on $A$, $B$ and $C$, we can appeal to  a version of the Center Manifold Theorem to aid in the proof. 
\smallskip

\begin{proof}To   adapt the necessary results from center manifold theory,  it is useful to formulate the problem in the state space $\sX = \sZ\times \sW$. Namely, with $u=\g(w)$ we have
\begin{align}
\label{eq8}
\dot X  &= \sA X + \sF(X) , \ X(0) = X_0,  \\[2ex]
 X &= \begin{pmatrix}z\\w\end{pmatrix},\ \ \sF(X)=\begin{pmatrix}f(z)+\wh{\g}(w)\\ 0 \end{pmatrix}, \nn \\[2ex]
 \sA &= \begin{pmatrix}A& (\Bin \G+\Bd P)\\0&S \end{pmatrix}, \nn \\[2ex]
 e& = C z - Q w .\nn
\end{align}

\smallskip
Under Assumption \ref{ass1} we have: for any fixed and continuously differentiable $\g$ and $f$, $\sF$ is continuously differentiable in $\sX^\a = \sZ^\a \times \sW$ for some $\a>0$, 
$$\sF(0)=0,\ \  \sF_X(0)=0$$
and the operator $(-\sA)$ is a sectorial operator, and therefore generates an analytic semigroup $\sT(t)$. Furthermore, 
$$\s(\sA) = \s(\sA)\cup \s(S), \ \ \s(\sA)\cap \s(S)=\emptyset ,$$
and for some $\b>0$, 
$$\s(\sA)\subset \bbc^-_{-\b}=\{\z: \re(\z)\leq -\b\}, \ \ \  \s(S)\subset\bbc_0 =i\bbr.$$
\smallskip

According to  \cite[Theorem 6.2.1]{henry}, for every $\g$ (as above), $\sX = \sX_1\otimes \sX_2$ where $\sX_j$ are $\sA$ invariant subspaces.  Where $X_1=\sP(\sX)(\sX) \simeq \sW$ (isomorphic to) and where $\sP$ is the projection onto the eigenspace spanned  by the finitely many eigenvalues of $\sA_1=\sA\big|_{\sX_1}$ (from \cite{IB} we recall that neutral stability of the exosystem implies $\s(\sA)\cap\bbc_0 =\s(S)$ consists of finitely many eigenvalues on the imaginary axis each having geometric multiplicity one) and $\sX_2 = \sP_2 \sX = (I-\sP)\sX$ corresponds to the spectrum of  $\sA_2=\sA\big|_{\sX_2}$.

In particular, according to  \cite[Theorem 6.2.1]{henry} 
 there is a $C^1$ local invariant manifold $\Sigma$ defined in a neighborhood  $U$   of the origin in $X^\a$ 
and a mapping $\pi\, :\, \sX_1\ra \sX$ so that, after identifying $\sW $ and $\sX_1$, we have
\begin{equation}
\label{eq9s1}
\Sigma = \left\{\begin{pmatrix} \pi(w) \\w\end{pmatrix}\, :\, w\in \sW  \right\}.
\end{equation}
Also $\sA_2$ generates an exponentially stable analytic semigroup $\sT_2$ in $\sX_2$.
\smallskip

Invariance of $\Sigma$ implies that for initial data $X(0)= (\pi(w),w)^\tra\in \Sigma$ the solution to \rr{eq8} satisfies $z(t)=\pi(w(t))$ for all $t\geq 0$.  In particular, for all time $t$,   on $\Sigma$ we have 
$$\begin{pmatrix}z(t)\\w(t) \end{pmatrix}   = \begin{pmatrix}\pi(w(t)) \\ w(t) \end{pmatrix} ,  $$
and therefore, from the equations of motion, we have
\begin{equation}
\label{invar1}
 \deriv{z}{t} = A\pi(w) +f(\pi(w))+B\g(w) + Pw ,
\end{equation}
and, by the chain rule, 
 \begin{equation}
\label{invar2}
\deriv{z}{t} = \deriv{\pi(w)}{t} = \pderiv{\pi}{w}\deriv{w}{t} = \pderiv{\pi}{w} S w. 
\end{equation}

Therefore, invariance of $\Sigma$ is equivalent to

\begin{equation}
\label{eq10s1}
  \pderiv{\pi}{w} S w=A\pi(w) +f(\pi(w))+B\g(w) + Pw. 
\end{equation}
\smallskip
Furthermore, for initial data $(z_0,w_0)^\tra \in U$, i.e., sufficiently small, the solution $(z,w)^\tra$ satisfies 
\begin{equation}
\label{eq11s1}
\|z(t) -\pi(w(t))\|_\a \leq K e^{-\a_0 t } \|z_0-\pi(w_0)\|_\a,
\end{equation}
 where $K$ depends on $\a$. That is, under our assumptions, the center manifold is  locally exponentially attractive.\\[2ex]
 
 \centerline{\includegraphics[width=2.95in]{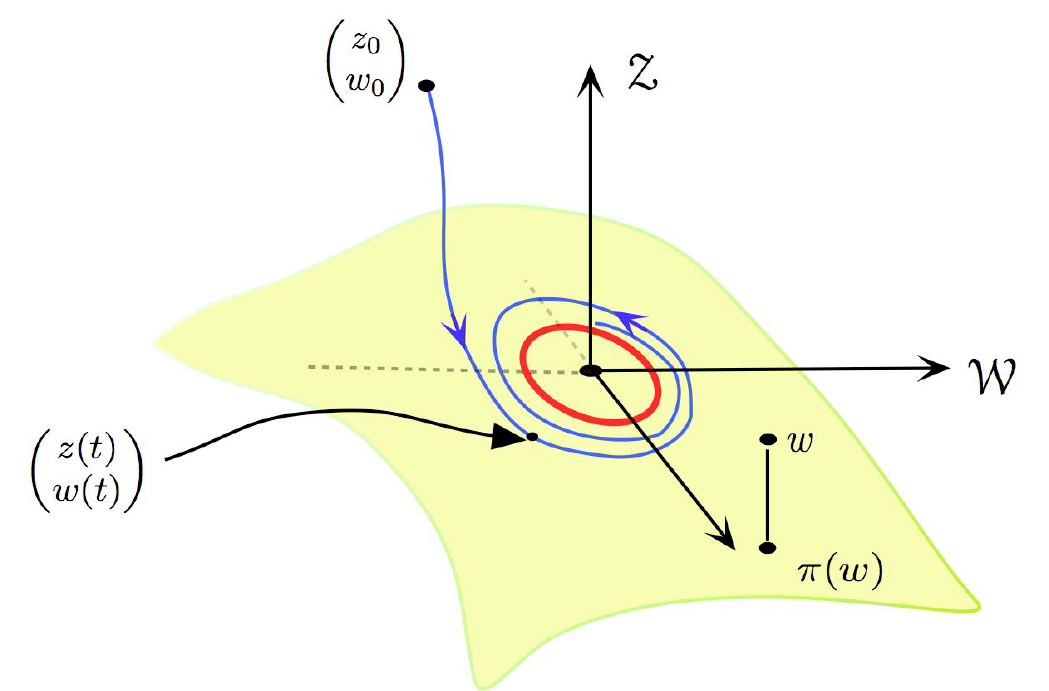}}
 
 \centerline{Figure 1: Center Manifold is Attractive }
\vskip.25in

\begin{remark}
In the general theory of dynamical systems,  center manifolds are not unique. However, as in \cite{IB}, an $\omega$-limit argument, using neutral stability  of the exosystem, implies uniqueness of  the center manifold $\Sigma$, given any choice of $\g$.
\end{remark}

By the same $\omega$-limit argument and the attractivity of the center manifold, it follows that   if  $u = \g(w)$ solves the output regulation problem, then since the center manifold is error zeroing, 
\begin{equation}
\label{eq12s1}
  e= Qw - C \pi(w)  = 0, 
\end{equation}
must hold. In particular, the solvability of the output regulation problem implies the solvability of the regulator equations \rr{sfrp1} and \rr{sfrp2}.

Conversely, if the regulator equations hold, then $z=\pi(w)$ is an invariant manifold which, by reversing the arguments above, must be the  center manifold $\Sigma$ for the closed-loop system. 
   Since $\Sigma$ is error zeroing and exponentially attractive, 
we have
\begin{align}
\|e(t)\|_\sY& = \|C z(t)- Q w(t)\|_\sY\nn =  \|[C z(t) -C \pi(w(t)) ]  
 +[C \pi(w(t))-Q w(t)]\|_\sY \nn \\[1ex]
&= \|C z(t) -C \pi(w(t) )\|_\sY   = \|C ( z(t)-\pi(w(t) )\|_\sY \nn \\[1ex]
&\leq c_\a\|z(t)-\pi(w(t))\|_\a \nn   \leq c_\a K e^{-\a_0 t } \|z_0-\pi(w_0)\|_\a \ra 0\ \ \text{ as }\  t\ra \i.\nn
\end{align}
 This  concludes the proof of Theorem \ref{thm1}. 
\end{proof}

\begin{remark} As first observed for the lumped case in \cite{IB}, the  dynamics on the invariant manifold $\sM$ is a copy, via $\pi$, of the dynamics of the exosystem and is therefore a model of the exosystem itself. This observation is the basis for the {\em internal model principle} for output regulation near an equilibrium, as developed in detail in  \cite{IsBook2}, \cite{BIbook} and more recently in \cite{Huang}. For lumped nonlinear systems, this has been extended to the nonequilibrium case in \cite{BI1}, \cite{BI2}, \cite{dutch}.
\end{remark}

\section{Dynamic Regulator Equations}\label{sec_rdreg}

In this section  we begin to develop a series of results, first introduced in \cite{aulisa2015practical}, providing  practical approximation  methods that extend the results from Section \ref{sec_regeq}  allowing for  analysis of very general regulation problems.

We want to return to the regulator equations in \rr{sfrp1} and \rr{sfrp2} and consider a dynamic version of these equations. 
In particular from  \rr{invar1} and \rr{invar2} we see that
\begin{equation}
\label{invar3}
 \deriv{\pi(w)}{t} = A\pi(w) +f(\pi(w))+\Bin \g(w) + \Bd Pw,
\end{equation}
and
 \begin{equation}
\label{invqr4}
C \pi(w) =Qw .
\end{equation} 
We now introduce the state variable $\ol{z}  = \pi(w)$ and  consider the particular case $\ol{z} (0) \eqdef \ol{z}_0 = \pi(w_0)$ where $w_0$ is the unique initial condition, guaranteed by the center manifold theorem,  so that $r = Qw$,  $d=Pw$.  Then let $\g$ be the  resulting control $u=\g(w)$, obtained by solving the regulator equations. Then, under the assumptions made in Section \ref{sec_regeq}, existence of $u=\g(w)$ is guaranteed by the center manifold theorem.

The result is a dynamical system, supplemented by the error zeroing  constraint given in \rr{dc3}.  This system  is completely equivalent to the regulator equations in \rr{sfrp1} and \rr{sfrp2} in the case when $r$ and $d$ are known. 
\begin{align} 
 \deriv{\ol{z}}{t} &= A\ol{z} +f(\ol{z})+\Bin \g  + \Bd d, \label{dc1}\\
 \ol{z}(0) & = \pi(w_0),  \label{dc2}\\
 C\ol{z}(t) & = r(t) \ \ \text{ for all } t\geq 0. \label{dc3}
\end{align}
Notice, in the system \rr{dc1}--\rr{dc3},  we do not know $\pi$ (without solving the regulator equations), which means that we have only rewritten the regulator equations. 

The regulator equations  \rr{sfrp1} and \rr{sfrp2} and, equivalently, \rr{dc1}--\rr{dc3} are usually very difficult to solve.  This fact is the primary motivation for seeking an approximation methodology that can deliver accurate tracking and disturbance rejection and is easily implemented in practice. 
In particular,  the dynamic form  \rr{dc1}--\rr{dc3} presents us with an opportunity to discover a means to obtain accurate approximate controls. This forms the main objective of the following work.  

In what follows, instead of considering \rr{dc1}--\rr{dc3}  we will replace the (unknown)  initial condition $ \pi(w_0) $ in \rr{dc2} by an initial condition $\ol{z}_0$ which is only assumed to be close to $ \pi(w_0)$. Notice that the geometric method is a local theory assumed to be studied in a neighborhood of an equilibrium  (as is the center manifold theorem). In this work we assume that equilibrium to be zero so we need only to choose a small enough initial condition. In this case,  we reformulate the control problem as follows:

Given a reference signal $r(t)$ and disturbance $d(t)$, find state variable $\ol{z}$ and control $\ol{u}$ so that, for any (sufficiently small) initial condition $\ol{z}_0$ the following system, which we refer to as the {\em dynamic regulator equations}, is satisfied 
\begin{align} 
 \deriv{\ol{z}}{t} &= A\ol{z} +f(\ol{z})+\Bin \ol{u}  + \Bd d, \label{dc4}\\
 \ol{z}(0) & = \ol{z}_0,  \label{dc5}\\
 r(t) &- C\ol{z}(t) \ra 0, \ \text{ as  }\ t\ra \i. \label{dc6}
\end{align}
As we have already seen in Section \ref{sec_regeq}, assuming the regulator problem is solvable, for any sufficiently small initial condition $\ol{z}_0$,  the system \rr{dc4}--\rr{dc6} has a solution  $z=\pi(w)$ and $\ol{u}$  given as solutions of the regulator equation \rr{sfrp1}--\rr{sfrp2}.   In this way, using the control $\ol{u}$ the system  \rr{dc4}--\rr{dc6}, will result in an error $e(t)$, that will  asymptotically approach $0$ as $t$ goes to infinity.

But, as will be seen later, an  attempt to solve \rr{dc4}--\rr{dc6} directly produces a singular system requiring a regularization.  Solving the regularized system for $\ol{z}$ and $\ol{u}$ produces an approximate control  so,  we cannot expect to satisfy \rr{dc6}.  Therefore, we obtain an error which we denote by $e_0(t) = r - C \ol{z} (t)$.  We expect this error to be small, but a proof that this is true would be required. 

The regularization will lead us to introduce a sequence of systems that will, with each step, produce increasingly more accurate controls. To begin the process let us replace $\ol{z}$ and $\ol{u} $ by $\ol{z}^0$ and $\ol{u}^0 $, respectively, in  \rr{dc4}--\rr{dc6} to solve our first approximation 
\begin{align} 
 \deriv{\ol{z}^0}{t} &= A\ol{z}^0 +f(\ol{z}^0)+\Bin \ol{u}^0  + \Bd d, \label{dc4_0}\\
 \ol{z}^0(0) & = \ol{z}_0,  \label{dc5_0}\\
 \mbox{with a desired objective }\nn \\
 C\ol{z}^0(t) & = r(t), \ \ \text{ for all } t\geq 0. \label{dc6_0}
\end{align}
As we will see in Section \ref{sec_rdregeq}, the singularity of the system \rr{dc4_0} -- \rr{dc6_0}, does not allow us to achieve \rr{dc6_0} exactly.

 \section{Regularized Dynamic Regulator Equations}\label{sec_rdregeq}

 In order to see how  the system \rr{dc4_0}--\rr{dc6_0} is singular,
 let us choose a particular $\ol{z}_0$,  and consider the problem of solving    \rr{dc4_0}--\rr{dc6_0}  for the control $\ol{u}^0(t)$ by rewriting  Eq.~\eqref{dc4_0} as
\begin{align}
& \ol{z}^0(t) = A ^{-1}\bigg[ \ol{z}^0_t(t) -  \Bd d(t)  -   f(\ol{z}^0)\bigg]    -   A ^{-1} \Bin \ol{u}^0(t). \label{eq5}
\end{align}
Then,  apply  $C$ to both sides and try to enforce  Eq. \eqref{dc6_0} to obtain 
\begin{align}
&r(t) = C(\ol{z}^0)(t) = CA ^{-1} \bigg[ \ol{z}^0_t(t)  -   f(\ol{z})  - \Bd d(t) \bigg] - C(A ^{-1}) \Bin \ol{u}^0 (t).\  \label{eq7} 
\end{align} 
Let us now denote  the transfer function of the linearization of the system,  at $s\in \bbc$ by  $G(s) = C(sI-A)^{-1}\Bin$. Then we   define $G$ to be $G(0)$, so that
\begin{equation}
\label{G}
G =   C(-A ^{-1}) \Bin
\end{equation} 
With this we can solve   Eq.~\eqref{eq7}  to obtain
  we have
  \begin{equation}
\label{eq9}
\ol{u}^0 = G^{-1} \bigg[ r - CA ^{-1} \bigg( \ol{z}^0_t(t)  -   f(\ol{z})  - \Bd d(t) \bigg) \bigg]. 
\end{equation}

To see how the system is singular we   substitute \rr{eq9} back into \rr{dc4_0} in an attempt to obtain a dynamical system that does not contain $\ol{u}^0$. This produces   
$$
  \ol{z}^0_t(t)   = A\ol{z}^0(t) +\big( f(\ol{z}^0) +   \Bd d(t)\big) + \Bin  G^{-1} \bigg\{ r(t) -  CA ^{-1} \ol{z}^0_t(t) +  CA^{-1}\big( f(\ol{z}^0) +   \Bd d(t)\big)\bigg\} .  
$$
Let us denote 
 \begin{equation}
\label{BinB}
B \defeq \Bin G^{-1}  = \Bin [C(-A ^{-1}) \Bin]^{-1}.
\end{equation} 

The advantage of $B$ over $\Bin$ is that 
$$ C(-A)^{-1}B = C(-A)^{-1}\Bin G^{-1} = 1.$$

This allows us to simplify some complicated expressions down the road. In particular,  we can write
\begin{equation}
(I+BCA ^{-1})\ol{z}^0_t(t) =  A\ol{z}^0(t) + B  r(t)  + (I+B CA^{-1}) \big( f(\ol{z}^0) +   \Bd d(t)\big),  \label{eq10}
\end{equation}
 which appears to be a dynamic equation not containing $\ol{u}^0$.  Unfortunately, as has been  shown in \cite{aulisa2018analysis,aulisa2015practical,aulisa2019analysis}  this  system is singular     
  since the operator $(I+BCA ^{-1})$ in front of the time derivative term in \rr{eq10} is not invertible. Indeed, $(-BCA ^{-1})$ is idempotent due to 
 $$ (-BCA ^{-1})^2 =  (-B CA ^{-1}) (-\Bin G^{-1} CA ^{-1})= (-B CA ^{-1}) $$  and 
therefore   $(I+BCA ^{-1})$  is also idempotent with a nontrivial null space.  

To deal with the  singular system we introduce  a regularization  by introducing a parameter  $0< \b<1$,   replacing  Eq. (\ref{eq9}) by 
\begin{equation}
\ol{u}^0 =   \bigg[ r - CA ^{-1} \bigg( (1-\b) \ol{z}^0_t(t)  -   f(\ol{z})  - \Bd d(t) \bigg) \bigg]. \label{eq11}
\end{equation}
 
We then repeat the above calculations    using  Eq.~\eqref{eq11} instead of \rr{eq9}  
to  arrive at 
\begin{equation}
  (I+(1-\b) BCA ^{-1}) \ol{z}^0_t(t) =  A\ol{z}^0(t) + B  r(t)  + (I+B CA^{-1}) \big( f(\ol{z}^0) +   \Bd d(t)\big).   \label{eq13} \\
\end{equation}
Using the condition $-CA^{-1}B =I$ from \rr{BinB} it follows that $(I+(1-\b) BCA ^{-1})$ is   invertible    with inverse
\begin{equation} 
\big[I+(1-\b) BCA ^{-1}\big]^{-1} = I - \frac{(1-\b)}{ \b}  (BCA ^{-1}).\label{eq14} 
\end{equation}
Applying    $(I+(1-\b)  BCA ^{-1} )^{-1}$ to both sides of equation \rr{eq13},  using  formula \rr{eq14}, and after a few simple   calculations,  \rr{eq13} can be written as
\begin{equation}
\ol{z}^{1} _{t}(t) = A_{\b}\ol {z}^{1}(t)  + (I+BCA^{-1})\bigg(f(\ol{z}^0) + \Bd d(t)\bigg) +\frac{1}{\b}B r (t). \label{reg18} 
\end{equation}
where 
\begin{equation}
A_\b =  \left( A- \z  BC \right), \ \ \z = \frac{(1-\b)}{ \b} . \label{eq17}
\end{equation}
In obtaining formula \rr{reg18}  and \rr{eq17} we have used the relations 
  \begin{align*}
&\left(I -\z  (BCA ^{-1}) \right)A =  \left( A-\z  BC \right) \equiv A_\b,  \\
&\left(I -\z  (BCA ^{-1}) \right)B =  \frac{1}{\b} B, \\
&\left(I -\z  (BCA ^{-1}) \right) (I+BCA^{-1})   = I+BCA^{-1}, 
\end{align*}
which are easily verified, as are the following
\begin{equation}
 C A_\b^{-1} = \b C A^{-1},       \ \ \ A_\b^{-1}B = \b\,   A^{-1}B, \ \ 
  C A_\b^{-1}B  =  - \b I , \label{cons2} 
\end{equation}
 which are used many times in what follows.

\begin{remark}\label{rem_Ab}{\rm 
Under our Assumptions \ref{ass1} it follows from Theorem 3.10.11 in  \cite{staffans}  that $A_\b$ is a sectorial operator with domain $\mathrm{D}(A_\b) = \mathrm{D}(A )$ since $P = - \z BC\in \sL(\sZ^{s_c}, \sZ^{-s_b})$. This result also shows that $A_\b$ generates an analytic semigroup $S_\b(t)$ on $\sZ^\g$ for $s_c-1\leq \g\leq 1-s_b$. In particular, under our assumptions of $s_b$ and $s_c$, this is true for $\g=0$ so $S_\b(t)$ defines an analytic semigroup in $\sZ$.    
 
Furthermore, by Proposition 5.4 in  \cite{MasterThesis},  $S_\b(t)$ is an exponentially stable semigroup for $\b$ sufficiently close to $1$.  In particular, with $\z =(1-\b)/\b$, for $\b$ sufficiently close to $1$ there is a constant $q$ so that 
$$\|PR(\l,A)\|_{-s_b} \leq q <1 \text{  for } \l \in \bbc^{+}= \{\l \in \bbc\,: \re (\l)>0 \}.$$
This follows from 
$$\|PR(\l,A)\|_{-s_b}  =  \z \|  B C \sa{-s_c} \sa{s_c} R(\l,A)\|_{-s_b} \leq   \z  \|B\|_{-s_b} \|C  \|_{-s_c}  \sup_{\l\in \bbc^+} \|(-A)^{s_c} R(\l,A)\| \defeq q . $$
Here we have
 $$\|(-A)^{ s_c} R(\l,A)\| = \|(-A)^{-1+ s_c} (-A) R(\l,A)\| \leq \|(-A)^{-1+ s_c}\| \| A R(\l,A)\|$$ 
and by our assumptions on $A$,  $\| A R(\l,A)\|$ is  bounded for $\l\in \bbc^+$. 
In this case  the growth bound for $A_\b$ satisfies  $-\o_\b<0$ and there is an $M_\b\geq 1$ so that 
\begin{equation}
\label{S_bdd}
 \|S_\b(t)\| = \|  e^{A_{\b}t}\| \leq M_{\b}  e^{-\o_{\b}t} .
\end{equation}
}
\end{remark}


At this point we make some statements about a slightly different approach taken in  \cite{aulisa2018analysis,aulisa2015practical,aulisa2019analysis}. The analysis in these works is completely equivalent to the  approach taken in  the above analysis.

\begin{remark}\label{rem_wtz}
In particular in \cite{aulisa2018analysis,aulisa2015practical,aulisa2019analysis}, 
the authors introduced a new variable 
$$\wt{z}^0 = A ^{-1} \bigg[ (1-\b) \ol{z}^0_t(t)  -   f(\ol{z}^0)  - \Bd d(t) \bigg]$$
which is the same as
\begin{equation}
\label{zt}
(1-\b) \ol{z}^0_t(t)  = A \wt{z}^0 +  f(\ol{z}^0) + \Bd d(t),
\end{equation}
 which    allowed them  to write 
\begin{equation}
\ol{u}^0(t) = G^{-1} \big[ r(t) -  C\wt{z}^0 \big].\   \label{eq9a}
\end{equation}
 
With this notation, in order to obtain the approximate control $\ol{u}^0$ in \rr{eq9}  they solved the  regularized  system  consisting of a pair of equations involving unknown variables  $\ol{z}^0$, $\wt{z}^0$ and $\ol{u}^0$, which satisfy 
 \begin{align}
& \ol{z}^0_{t}(t) = A \ol{z}^0 (t)  + f(\ol{z}^0)  + \Bd d(t) + \Bin  \ol{u}^0 (t), \label{e1}   \\
& (1-\b) \ol{z}^0 _{t}(t) = A \wt{z}^0(t)  + f(\ol{z}^0)  +\Bd d(t), \label{e2} \\
& \ol{z}^0 (0)	    = \ol{z}^0_0. \label{e3} 
\end{align}
The result then gives 
\begin{equation}
\label{olg1}
\ol{u}^0 = G^{-1}\big( r(t) - C \wt{z}^0(t)  \big).
\end{equation}

The system \rr{e1}--\rr{olg1} is easily solvable using off-the-shelf software. For example, the finite element software package Comsol, \cite{COMSOL}, can be used to solve this system for a wide range of practical applications when $A$ is given as an ordinary or partial differential operator.
\end{remark}


\section{Regularized Dynamic Controller}\label{sec_rdc}

Our objective now is to obtain a more standard form of our control design  as a dynamic controller. Some of material below is either taken from or consists of a modified version of the material found in \cite{aulisa2018analysis,aulisa2015practical,aulisa2019analysis}.  But the notion of dynamic controller as it is presented here is new. 
   
We have demonstrated that the regularized version of  \rr{dc4_0}-\rr{dc6_0} can be written as the dynamical system
\begin{align}
\ol{z}^{0} _{t}(t) &= A_{\b}\ol {z}^0(t) + I_0 \bigg( f(\ol{z}^0) +\Bd d(t) \bigg) +\frac{1}{\b}B r (t) ,  \label{rdc1} \\
\ol{z}^{0} (0)& =\ol{z}_0, \label{rdc2} 
\end{align}
where
\begin{equation}
\label{eq_d}
I_0  = (I+BCA^{-1}) .
\end{equation}
Solving the system \rr{rdc1}--\rr{rdc2} we then obtain a formula for the control (given in \rr{eq11}).
$$\ol{u}^0 = G^{-1} \bigg[ r - CA ^{-1} \bigg( (1-\b) \ol{z}^0_t(t)  -   f(\ol{z})  - \Bd d(t) \bigg) \bigg]. $$

A practical method  for approximating $\ol{u}^0$   is described above in Remark  \ref{rem_wtz}. 

 We reiterate an important point. Due to the introduction of the regularization $(1-\b)$  in Eq.~\eqref{eq11}  we no longer have 
 $r(t) - C(\ol{z}^0) (t)=0 $ for all $t$.   We denote the resulting difference as an error
 \begin{equation} 
e_0(t) = r(t) - C(\ol{z}^0) (t)  . \label{eqE1} 
\end{equation}

  We now proceed to use the above information to  simplify our formula for $\ol{u}^0$ given in   \rr{eq11}. In particular, using \rr{eq11} and the new representation for $\ol{z}^{0} _{t}(t)$ given in \rr{rdc1} we have
 \begin{align}
&\ol{u}^0(t)  =  G^{-1}\bigg[r(t) - CA ^{-1}\bigg( (1-\b)  \ol{z}^0_t(t)  - f(\ol{z}^0(t))  - \Bd d(t) \bigg) \bigg]  \label{eq12}\\
& = G^{-1}\bigg[ r(t) -  (1-\b)  CA ^{-1}\bigg\{ A_{\b}\ol {z}^0(t) + I_0\big(f(\ol{z}^0(t))  +  \Bd d(t)\big) +\frac{1}{\b}BG^{-1} r (t) \bigg\} \nn \\
&  \hspace*{.25in}    -  CA^{-1} \big(  f(\ol{z}^0(t))  + \Bd d(t)\big)  \bigg] \nn  \\
&=   G^{-1}\bigg[\left( 1 - \frac{(1-\b)}{\b}  CA ^{-1}BG^{-1} \right) r (t) -  (1-\b)  CA ^{-1} A_{\b}\ol {z}^0(t)  \nn \\
& \hspace*{.25in}    -  (1-\b)  CA ^{-1}I_0\big(f(\ol{z}^0(t))  +  \Bd d(t) \big)    +  CA ^{-1}\big( f(\ol{z}^0(t))  +  \Bd d(t) \big)  \bigg] \nn  \\
& = G^{-1}\bigg[ \frac{1}{\beta} r(t) -  \frac{(1-\b)}{\beta} C \ol{z}^0(t) + CA^{-1}\bigg(f(\ol{z}^0(t) ) +    \Bd d(t)\bigg)\bigg]. \label{eq12a}
\end{align}
Here we have used 
$$   CA^{-1}I_0 =CA^{-1}(I+BG^{-1}C(A)^{-1})   =0, \ \ CA^{-1}A_\b = \frac{1}{\b} C.$$  

Thus we have obtained a simple formula for the desired approximate control $\ol{u}^0$, that results in the error $e_0(t)$, when trying to track $r(t)$  
\begin{equation}
\label{gam0}
\ol{u}^0(t)  = G^{-1}\bigg[ \frac{1}{\beta} r(t) -  \frac{(1-\b)}{\beta} C \ol{z}^0(t) + CA^{-1}\bigg(f(\ol{z}^0(t) ) +    \Bd d(t)\bigg)\bigg].
\end{equation}

Combining these results we obtain the {\em Regularized Dynamic Controller} for first iteration
\begin{align}
\ol{z}^{0} _{t}(t) &= A_{\b}\ol {z}^0(t) + I_0 \bigg( f(\ol{z}^0) +\Bd d(t) \bigg) +\frac{1}{\b}B r (t) ,  \label{rdc3} \\
\ol{z}^{0} (0)& =\ol{z}_0, \label{rdc4} \\
\ol{u}^0 &= G^{-1}\bigg[ \frac{1}{\beta} r(t) -  \frac{(1-\b)}{\beta} C \ol{z}^0(t) + CA^{-1}\bigg(f(\ol{z}^0(t) ) +    \Bd d(t)\bigg)\bigg]. \label{rdc5}
\end{align}

\begin{remark}\label{rem_z_0} 
A substantial difference between what we have done here and what was done in \cite{aulisa2018analysis,aulisa2015practical,aulisa2019analysis} is that  we have allowed a more arbitrary initial condition $\ol{z}_0$ in \rr{dc5} and subsequently in \rr{rdc4}. In the previous works the authors chose very specific initial conditions. Different initial data result in different transient behavior. In the error analysis, this yields more complicated formulas, due to extra terms that decay exponentially to zero, but do not otherwise change the asymptotic behavior. The reason we have done things differently here is that for the delay equations the procedure used in \cite{aulisa2018analysis,aulisa2015practical,aulisa2019analysis} cannot be applied straightforwardly.  

In order for the reader to get a better understanding of the difference, we give a very brief outline of how the specific initial condition is chosen in the earlier works.  

In these other works the authors obtain a specific initial condition by first solving  a set-point regulation problem to track the time independent  reference signal, $r(0)$, and reject the time independent  disturbance, $d(0)$.  
Namely,  they solve the following stationary problem to obtain $\ol{u}^0$ and $\ol{z}^{0}$
\begin{align}
& 0 =  A \ol{z}^{0} + B\ol{u}^0  +  \Bd d(0) + f(\ol{z}^0),  \label{eqz0} \\
& r(0) = C( \ol{z}^{0} ). \label{eqz0x} 
\end{align}  
Solving the above for $\ol{z}^0$ and $\ol{u}^0$  is very easy to carry out numerically using, for example, an off the shelf finite element software package like Comsol \cite{COMSOL}. With this approach the first iterative error (using a capital $E$ for error to avoid confusion with the $e$ used in this work)
\begin{equation}
E_0(t) = r(t)-r(0) .\label{E0} 
 \end{equation}
 and similarly the perturbation with respect to its initial value for the disturbance $d$
 \begin{equation}
\label{wtd}
D_H(t)=d(t)-d(0).
\end{equation}
The solution  $ \ol{z}^{0}$ is  used as an accurate approximation for the unknown  initial condition $ \pi(w_0)$ in \rr{dc2}.
Another advantage of this approach is that $E_0(0)=0$ which, once again, simplifies some of the later computations. 
In particular, in carrying out the estimation of the errors in Section \ref{sec_error} below. At each step we do not have $e_j(0)=0 $ as it is in the case in which this special initial initial condition is used.  Namely the analogs $E_j(t)$ of $e_j(t)$ in \cite{aulisa2018analysis,aulisa2015practical,aulisa2019analysis} all satisfy $E_j(0)=0$. 

Beginning with the set-point problem also introduces another difference in the iterative scheme. Namely, in the iterative scheme introduced in Section \ref{iterate} below, the first terms in $\ol{z}_n$ and $\ol{u}_n$   
$$\ol{z}_{n} = \sum_{j=0}^n \ol{z}^j, \ \ \ \ol{u}_{n} = \sum_{j=0}^n\ol{u}^{j} $$
are now already fixed.  This  translates into the equations generating  $\ol{z}^1$ and $\ol{u}^1$ become essentially the same as the equations for $\ol{z}^0$ and $\ol{u}^0$ in the approach used in this work.   The equations generating $\ol{z}^j$ and $\ol{u}^j$ for $j>1$ are the same for both approaches. 

The result of allowing for a more general initial condition to start with  is that one needs to deal  with more complicated transients that decay rapidly to zero at each iteration.  
 
\end{remark}

\section{An Iterative Dynamic Controller}\label{iterate}

  We now define an  iteration procedure that allows us to obtain increasingly more accurate approximate controls. We begin by defining a sequence  of state variables and controls, beginning with  $\ol{z}_0 = \ol{z}^0$, $\ol{u}_0 = \ol{u}^0$ from \rr{rdc3} -- \rr{rdc5}, and for $n\in \bbz^+$, with $n\geq1$
 \begin{equation}
\label{ziter_gamiter_n}
\ol{z}_{n} = \sum_{j=0}^n \ol{z}^j, \ \ \ \ol{u}_{n} = \sum_{j=0}^n\ol{u}^{j} 
\end{equation}

 For a fixed $n$,  we define a sequence  of  {\em Regularized Systems}   for $1\leq j\leq n$
   \begin{equation}
\deriv{}{t}\ol{z}^j(t) =A_\b \ol{z}^j(t) + L_{0}    F_j    + \frac{1}{\beta}B    e_{j-1}(t), \ \ \ol{z}^j(0)=0  ,  \label{rdcj}
\end{equation}
  where the nonlinear terms are given by
  \begin{equation}
\label{Fj}
F_j = f(\ol{z}_{j}(t)) -  f(\ol{z}_{(j-1)}(t)),
\end{equation}
and the sequence of errors $e_j(t)$ are  given by
  \begin{equation}
\label{eq_ej}
e_j(t) = e_{j-1}(t)-C \ol{z}^j.
\end{equation}

These systems produce  a sequence of  controls
   \begin{equation}
\ol{u}^j(t)  = G^{-1} \bigg[ \frac{1}{\beta} e_{j-1}(t) -  \frac{(1-\beta)}{\beta}C  \ol{z}^j(t)  +   CA^{-1}  F_j  \bigg].  \label{eq:iter_j_control}
\end{equation}

Solving the systems for $j=1, \cdots , n$ we obtain a new control  $\ol{u}_{ n} $, as described in \rr{ziter_gamiter_n},  which delivers an error $e_n(t)$.

In order to see where these systems come from let us consider the case $n=1$, so that $\ol{z}_{1}= \ol{z}^0+\ol{z}^1$ and  $\ol{u}_{1}= \ol{u}^0+\ol{u}^1$. We assume that  $\ol{z}^0$ and $\ol{u}^0$   are solutions of
\rr{dc4_0}--\rr{dc6_0} and $\ol{z}_{1}$ and $\ol{u}_{1}$ are solutions of \rr{dc4}--\rr{dc6}. Then let us consider  what system $\ol{z}^1$ and $\ol{u}^1$ satisfy.  To this end we substitute $\ol{z}_1$ and $\ol{u}_1$ into \rr{dc4}:
$$\deriv{}{t}  \ol{z}_1  = A\ol{z}_1 +f(\ol{z}_1)+\Bin \ol{u}_1  + \Bd d, $$
which implies
$$\deriv{}{t} (\ol{z}^0+\ol{z}^1 )  = A(\ol{z}^0+\ol{z}^1 )  +f(\ol{z}^0+\ol{z}^1 )+\Bin(\ol{u}^0+\ol{u}^1)  + \Bd d ,$$
and since $\ol{z}^0$ and $\ol{u}^0$ satisfy \rr{dc4_0}
$$  \deriv{\ol{z}^0}{t}  = A\ol{z}^0 +f(\ol{z}^0)+\Bin \ol{u}^0  + \Bd d.$$
we see that  $\ol{z}^1$ and $\ol{u}^1$ satisfy
$$\deriv{}{t}  \ol{z}^1    = A \ol{z}^1   +f(\ol{z}^0+\ol{z}^1 ) -f(\ol{z}^0) +\Bin \ol{u}^1  $$
which, due to \rr{Fj}, can be written  as
$$\deriv{}{t}  \ol{z}^1    = A \ol{z}^1   +F_1 +\Bin \ol{u}^1 . $$
 Notice, in particular, that the disturbance $d$ no longer appears.
 
 Next we consider \rr{dc5} and \rr{dc5_0} which give
 $$\ol{z}_1(0)  = \ol{z}^0(0)+\ol{z}^1(0) = \ol{z}_0  \text{  and  }\ol{z}^0(0)   = \ol{z}_0 ,  $$
 which implies
 $$\ol{z}^1(0)  = 0.$$
 Finally, recalling our desired goal in  \rr{dc6} and the fact that $r-C \ol{z}^0 = e_0$ we want
 $$  C\ol{z}(t)  =  C ( \ol{z}^0+\ol{z}^1)    = r(t). $$
 This would mean that 
 $$C \ol{z}^1     = r(t) - C  \ol{z}^0 = e_0(t).$$
 So we see that the tracking objective for the control $\ol{u}^1$  for the system
 \begin{align}
 \deriv{}{t}  \ol{z}^1    &= A \ol{z}^1   + F_1 +\Bin \ol{u}^1 \label{sys1_a}\\
\ol{z}^1(0)  &= 0 \label{sys1_b}
\end{align}
is
\begin{equation}
\label{sys1_c}
C \ol{z}^1 (t)    =  e_0(t).
\end{equation}

Notice that the system \rr{sys1_a}--\rr{sys1_c} is identical (at least  in form) to \rr{dc4_0}--\rr{dc6_0} with the following differences:
\begin{enumerate}
\item There is a slightly different nonlinear term
$$ F_1 = f(\ol{z}^0+\ol{z}^1 ) -f(\ol{z}^0) \text{ in place of }   f(\ol{z}^0);$$
\item A different initial condition 
$$\ol{z}^1(0)  = 0  \text{ in place of }   \ol z^0(0) =  \ol{z}_0;$$
\item A different target reference signal
$$ C \ol{z}^1 (t)    =  e_0(t) \text{ in place of }   C(\ol{z}^0) (t) = r(t) ;$$
\item For the new system there is no disturbance term, i.e., 
$ d(t)=0$;
\item A new  desired control 
$$  \ol{u}^1  \text{ in place of }   \ol{u}^0 . $$
\end{enumerate}

As we have already seen, in attempting to solve for $\ol{u}^1$ in a system in the form \rr{sys1_a}--\rr{sys1_b} one encounters  a singular system. Therefore we  repeat all the calculations obtained  in our attempt to solve for $\ol{u}^0 $ (including the regularization involving $\b$  in \rr{eq11},  that lead to \rr{rdc1}--\rr{rdc2}). In this way we obtain a new regularized version of system for $ \ol{z}^1$, $ \ol{u}^1$ given by
\begin{align}
\ol{z}^{1} _{t}(t) &= A_{\b}\ol {z}^1(t) + I_0 F_1 +\frac{1}{\b}B e_0 (t) ,  \label{rdc1_1} \\
\ol{z}^{1} (0)& = 0. \label{rdc2_1} 
\end{align}

Once again we can repeat the steps \rr{eq12}--\rr{eq12a} which produced  the control $\ol{u}^0$ given in \rr{gam0}. Carrying out these same calculations  gives us the following  formula for $\ol{u}^1$.
\begin{equation}
\label{gam1}
\ol{u}^1(t)  = G^{-1}\bigg[ \frac{1}{\beta} e_0(t) -  \frac{(1-\b)}{\beta} C \ol{z}^1(t) + CA^{-1}F_1 \bigg].
\end{equation}
which is exactly the formula given in \rr{eq:iter_j_control} in the case $j=1$. 

It should be clear from the above discussion how the remaining systems for $j>1$ in \rr{rdcj}
 and the resulting errors $e_j(t)$ in \rr{eq_ej} and the controls $\ol{u}^j(t) $ in \rr{eq:iter_j_control} arise.

The main difficulty at this point is that we have not provided any evidence to support the claim that as we increase the number of iterations the resulting error $e_n(t)$ will decrease.  This was established for many examples in a series of articles and the book \cite{aulisa2015practical}. A fairly complete discussion for the analysis of the error in the linear case was provided in  the paper \cite{aulisa2019analysis}.   The nonlinear case is much more involved  but some general  results are obtained in \cite{aulisa2018analysis} for the first three iterations. In the next section we outline the basic results obtained in the these works.

Before beginning to examine the estimates of errors, we conclude this section with a discussion of the iterative method in the case of the notation introduced in the Remark \ref{rem_wtz}. These are essentially the iterative schemes studied in \cite{aulisa2015practical,aulisa2019analysis,aulisa2018analysis}. 

\begin{remark}\label{rem_iter_wtz}
 For $n\geq 1$ we introduce the same functions $\ol{z}_{n}$ and $\ol{u}_{n}$  given in  \rr{ziter_gamiter_n}. But using the notation for the systems given in \rr{e1}--\rr{e3} we obtain the analogs of 

 \begin{align}
& \ol{z}^j_{t}(t) = A \ol{z}^j (t)  + F_j    + \Bin  \ol{u}^j (t), \label{e1_iter}   \\
& (1-\b) \ol{z}^j _{t}(t) = A \wt{z}^0(t)  + F_j ,   \label{e2_iter} \\
& \ol{z}^j (0)	    = 0, \label{e3_iter} 
\end{align}
where
$$\wt{z}^j = A ^{-1} \bigg[ (1-\b) \ol{z}^j_t(t)  -  F_j   \bigg].$$

in this case we still have the same errors 
 \begin{equation} 
E_j(t) = E_{j-1}(t) - C(\ol{z}^j) (t)  , \label{eqEj} 
\end{equation}
 and the above system  delivers a new approximate control 
\begin{equation}
\label{olg1_iter}
\ol{u}^j= G^{-1}\big( E_{j-1} - C \wt{z}^j(t)  \big).
\end{equation}

\end{remark}

\section{Error Estimates}\label{sec_error}

Using the systems described in Remarks \ref{rem_wtz} and \ref{rem_iter_wtz}, 
the paper \cite{aulisa2019analysis} contains a complete analysis of error estimates for the linear regulation problem (when $f=0$) with tracking and disturbance rejection and the paper  \cite{aulisa2018analysis} provides some partial results for the nonlinear problem for tracking but without disturbances.  

The analysis of errors for the iterative scheme is rather lengthy and technical.  But the main calculations that allow for obtaining estimates of the errors derive from the dynamical systems in 
\rr{rdc1}--\rr{rdc2} and \rr{rdcj}.  Even though the  analysis contained in \cite{aulisa2019analysis} is quite involved  we will provide an outline in Subsection \ref{subsec_linear} for some of the main ingredients  and allow the reader to look at the full  paper to fill in all the details.  The nonlinear case is much more involved and while only partial results were obtained in  \cite{aulisa2018analysis}, these results provide very good estimates of the errors for the first three iterations. In practice, no more than one or two iterations are typically necessary to achieve satisfactory results.

In either the linear or nonlinear case the starting point is  applying the variation of parameters formula. Let us begin with  \rr{rdc1}--\rr{rdc2} from which we obtain
\begin{equation}
\label{gen3}
 \ol{z}^{0}(t)=e^{A_\b t} \ol{z}_0 +  \int_0^t e^{A_\b (t-\tau)} \left( \frac{1}{\b} B r(\tau) + I_0    \Bd d(\tau) \right)  \, \dd\tau +   \int_0^t e^{A_\b (t-\tau)}  I_0   f( \ol{z}^{0}(\tau)   \, \dd\tau .
 \end{equation} 
Next we integrate by parts in the first integral only  
 \begin{align*}
&\ol{z}^{0}(t)  = e^{A_\b t} \ol{z}_0 + \int_0^t e^{A_\b (t-\tau)} \bigg( \frac{1}{\b} B r(\tau) + I_0 \Bd d(\tau) \bigg)\, \dd\tau +   \int_0^t e^{A_\b (t-\tau)}  I_0    f( \ol{z}^{0}(\tau) ) \, \dd\tau \\
    &= e^{A_\b t} \ol{z}_0 + \int_0^t \deriv{}{\tau}\bigg[(-A_\b)^{-1}e^{A_\b (t-\tau)}\bigg] \bigg( \frac{1}{\b} B r(\tau) + I_0 \Bd d(\tau) \bigg)\, \dd\tau \\
    &  \hspace*{.25in} +   \int_0^t e^{A_\b (t-\tau)}  I_0    f( \ol{z}^{0}(\tau) ) \, \dd\tau   \\
    &  = e^{A_\b t} \ol{z}_0 +\bigg[(-A_\b)^{-1}e^{A_\b (t-\tau)}\bigg] \bigg( \frac{1}{\b} B r(\tau) + I_0 \Bd d(\tau) \bigg)\bigg|_0^t \\
    & \hspace*{.25in} - \int_0^t \bigg[(-A)^{-1}e^{A_\b (t-\tau)}\bigg] \deriv{}{\tau}\bigg( \frac{1}{\b} B r(\tau) + I_0 \Bd d(\tau) \bigg)\, \dd\tau +   \int_0^t e^{A_\b (t-\tau)}  I_0    f( \ol{z}^{0}(\tau) ) \, \dd\tau   \\
    & = e^{A_\b t} \ol{z}_0 + (-A_\b)^{-1} \bigg[  \frac{1}{\b}B r(t) +  I_0 \Bd d(t) \bigg] + A_\b^{-1}e^{A_\b t} \bigg[  \frac{1}{\b}B r(0) +  I_0 \Bd d(0) \bigg]  \\
    & + \int_0^t  A_\b^{-1}e^{A_\b (t-\tau)} \bigg( \frac{1}{\b} B r'(\tau) + I_0 \Bd d'(\tau) \bigg)\, \dd\tau +   \int_0^t e^{A_\b (t-\tau)}  I_0    f( \ol{z}^{0}(\tau) ) \, \dd\tau .
\end{align*}

Next we apply $C$ to both sides and use   three easily verified  identities:
$$ C(-A)^{-1}B =1, \ \  \frac{1}{\b} CA_\b^{-1} = CA^{-1}, \ \  C(-A)^{-1}I_0 = C(-A)^{-1} + C(-A)^{-1}BCA^{-1} = 0,$$
allowing us to write
\begin{align}
C\ol{z}^0 = & C e^{A_\b t} \ol{z}_0 +  r  + C A_\b^{-1}e^{A_\b t} \bigg[  \frac{1}{\b}B r(0) +  I_0 \Bd d(0) \bigg] \label{eeq1} \\
  & + C \int_0^t  A_\b^{-1}e^{A_\b (t-\tau)} \bigg( \frac{1}{\b} B r'(\tau) + I_0 \Bd d'(\tau) \bigg)\, \dd\tau +   C\int_0^t e^{A_\b (t-\tau)}  I_0    f( \ol{z}^{0}(\tau) ) \, \dd\tau .  \nn
\end{align}

Rewriting \rr{eeq1} and recalling that $ e_0 = r-C\ol{z}^0$ we have
\begin{align}
 &e_0(t) =  C e^{A_\b t} \ol{z}_0 -C A_\b^{-1}e^{A_\b t} \bigg[  \frac{1}{\b}B r(0) +  I_0 \Bd d(0) \bigg]  \\
 &- C \int_0^t  A_\b^{-1}e^{A_\b (t-\tau)} \bigg( \frac{1}{\b} B r'(\tau)  
  + I_0 \Bd d'(\tau) \bigg)\, \dd\tau  
     -  C \int_0^t e^{A_\b (t-\tau)}  I_0    f( \ol{z}^{0}(\tau) ) \, \dd\tau . \label{eeq2}
\end{align}

To simplify  many of the estimates obtained later we define the following class of functions.
\begin{definition}\label{def5}
We denote by $\calg$ all functions   $g(t)$ in the form $p(t)e^{-\o_\b t}$ 
where $p(t)$ is a polynomial in $t$.  Notice that functions in the set $\calg$ 
go to zero exponentially fast as $t$ goes to infinity.
\end{definition}
In particular, we see that
\begin{equation} G_0(t) \defeq  C e^{A_\b t} \ol{z}_0  -C A_\b^{-1}e^{A_\b t} \bigg[  \frac{1}{\b}B r(0) +  I_0 \Bd d(0) \bigg] \label{G0eq} \end{equation}
 decays to $0$ exponentially as $t$ goes to $\i$ and is clearly bounded above by a function $g_0 \in \calg$. 
 
We comment that with a more auspicious choice of initial condition (as discussed in Remark \ref{rem_z_0}) and a slightly different beginning to our iterative method we could have avoided having the term $G_0(t)$. 
  
Returning to our iterative analysis in Section \ref{iterate}, we have $\ol{z}_{n}$, $\ol{u}_{n}$ defined in \rr{ziter_gamiter_n} in terms of $\ol{z}^{j}$, $\ol{u}^{j}$ for $j=0, 1, \dots ,n$ and  $e_j(t) = e_{j-1}(t)-C \ol{z}^j$  given in \rr{eq_ej}.  
Recall that $\ol{z}^{j}$, $\ol{u}^{j}$ satisfy \rr{rdcj} 
$$\deriv{}{t}\ol{z}^j(t) =A_\b \ol{z}^j(t) + L_{0}    F_j    + \frac{1}{\beta}B    e_{j-1}(t), \ \ \ol{z}^j(0)=0 $$
and from \rr{eq:iter_j_control}
$$\ol{u}^j(t)  = \bigg[ \frac{1}{\beta} e_{j-1}(t) -  \frac{(1-\beta)}{\beta}C  \ol{z}^j(t)  +   CA^{-1}  F_j  \bigg]. $$

Applying  the variation of parameters formula to \rr{rdcj}  we obtain
\begin{equation}
\label{genj}
 \ol{z}^{j}(t)=  \int_0^t e^{A_\b (t-\tau)}  \frac{1}{\b} B e_{j-1}(\tau)\, \dd\tau +   \int_0^t e^{A_\b (t-\tau)}  I_0   F_j(\tau)  \, \dd\tau .
 \end{equation} 
Next we repeat the integration by parts carried out above to obtain 
\begin{align}
C\ol{z}^j = &   e_{j-1}(t)    + C A_\b^{-1}e^{A_\b t} \bigg[  \frac{1}{\b}B e_{j-1} (0)   \bigg] \label{eeqj} \\
  & + C \int_0^t  A_\b^{-1}e^{A_\b (t-\tau)} \bigg( \frac{1}{\b} B e_{j-1}'(\tau)  \bigg)\, \dd\tau +   C \int_0^t e^{A_\b (t-\tau)}  I_0   F_j (\tau)  \, \dd\tau .  \nn
\end{align}
Rewriting \rr{eeqj} and recalling that $ e_{j} = e_{j-1}-C\ol{z}^j$ we have
\begin{align}
e_j(t)  &  = - C A_\b^{-1}e^{A_\b t} \bigg[  \frac{1}{\b}B e_{j-1} (0)   \bigg]   \\
    & -  \int_0^t  A_\b^{-1}e^{A_\b (t-\tau)} \bigg( \frac{1}{\b} B e_{j-1}'(\tau)  \bigg)\, \dd\tau -   \int_0^t e^{A_\b (t-\tau)}  I_0   F_j (\tau)  \, \dd\tau . 
\end{align}

In order to simplify the  outline of the error analysis let us define the following operators
\begin{align}
 K(t)=& - \frac{1}{\b} C A_\b^{-1} e^{A_\b t} B,\ \ K_d(t)  = -CA^{-1}_{\beta} e^{A_\beta t}I_0 \Bd,  \label{eq29}\\
   H =& -  C   e^{A_\b t}I_0 , \ \ I_0 = (I+BCA^{-1}). \label{eq29_b}
\end{align}
We also note that once again the norm of the operator
$$G_j(t) = - C A_\b^{-1}e^{A_\b t} \bigg[  \frac{1}{\b}B e_{j-1} (0)   \bigg] $$ 
is bounded above by a function $g_j\in \calg$. 
These operators will be useful in our analysis of the errors $e_n$ in our iterative scheme.

In particular, for $n\geq 1$ we can write 
\begin{align}
e_n(t)  &  = G_n(t)     +   \int_0^t  K(t-\tau)  e_{n-1}'(\tau)   \, \dd\tau +  \int_0^t H  (t-\tau)     F_n (\tau)  \, \dd\tau \nn \\
& = G_n(t)     + (K\ast e_{n-1}')(t) + (H\ast F_n)(t), \label{eq2_b}
\end{align}
where we have introduced the convolution operator defined for two functions $f$ and $g$   by
$$(f\ast g)(t) = \int_0^t f(t-\tau) g(\tau)\, d\tau.$$

With this notation, and recalling \rr{G0eq}
we can write \rr{eeq2} as
\begin{equation}
\label{eq2_a}
e_0(t) = G_0(t) + (K\ast r')(t) + (K_d\ast d')(t) + (H\ast F_0)(t).
\end{equation}

We now introduce the main results from two different works. The paper \cite{aulisa2019analysis} is concerned with the linear case $f(\cdot)=0$. The other paper \cite{aulisa2018analysis}, is concerned with the nonlinear case.  For linear systems there is a fairly complete discussion for estimating  the  errors $e_n(t)$. In the nonlinear case things are much more complicated.  The authors were able to obtain estimates for the errors $e_1(t)$, $e_2(t)$ and $e_3(t)$. In the nonlinear paper, only tracking was considered so there was no disturbance. In the presentation below for nonlinear systems we will present a brief discussion of  the estimates including disturbances. 

As a final remark in this section, in  the linear or even nonlinear case,   the iterative scheme delivers very accurate tracking and disturbance rejection after only one or two iterations. The size of the errors at each step are influenced by the choice of regularization parameter $\b$. Generally, a smaller $\b$ delivers smaller values of the $e_n$, but as $\b$ becomes smaller,  the operator $A_\b$ eventually becomes unstable, as discussed, for example, in \cite{aulisa2019analysis}.   Nevertheless, seldom are higher order iterations needed.

\subsection{Linear Systems}\label{subsec_linear}

The material in this section is a slightly modified version of the results in  \cite{aulisa2019analysis}.  
In this section we consider \rr{eq2_a} and \rr{eq2_b} in the case $f(\cdot)=0$.  In this case these equations simplify to
\begin{equation}
\label{eq2_a_fe0}
e_0(t) = G_0(t) + (K\ast r')(t) + (K_d\ast d')(t) ,
\end{equation}
\begin{align}
e_j(t)  &  =   G_j(t)     + (K\ast e_{j-1}')(t) . \label{eq2_b_fe0}
\end{align}

 To describe what these have to say in the present case we define some important values. 
We define
  \begin{equation}
\label{Dform}
D=   \int_0^\i \|  K(t) \|\, \dd t   , \ \ \ \ \ 
D_d=  \int_0^\i \|  K_d(t) \|\, \dd t  
\end{equation}
and set
$$ \mathfrak{D} = \frac{D_d}{D}. $$

The expressions $D$ and $D_d$ are finite as shown  in Remark 1.6 in \cite{aulisa2019analysis}. The following result  is  Theorem 1.1  in \cite{aulisa2019analysis}.

\begin{theorem}\label{main_thm}
 Let
 \begin{equation}
\label{eq_Cn}
 C_{n} =     \Big(  \max_{0\leq j \leq n} \limsup_{[0,\i)}\|r^{(j)}(t)\|_{\sY} +  \mathfrak{D}  \max_{0\leq j \leq n}  \limsup_{[0,\i)}\|d^{(j)}(t)(t) \|_{\sD} \Big)  . 
\end{equation}
 Then 
  \begin{equation}
\label{mainest}
\limsup_{t>0}  \|e_n(t)\|_\sY =  \limsup_{t>0}  \|r(t) -C(\ol{z}_n)(t)\|_\sY= \limsup_{t>0}  \|e_n(t)\|_\sY   \leq   D^n \  C_{n} . 
\end{equation}

Moreover, if there exist constants $\ol{\a}$ and $\ol{A}$ so that
$$C_n \leq \ol{A} \, \ol{\a}^n \  \text{    and    } \  D \ol{\a}  <1, $$
then
  $$ \limsup_{t>0}  \|E_n(t)\|_\sY  \leq  \ol{A}  (D\ol{\a})^n   \xrightarrow{n\ra \i} 0.$$
 In this case for sufficiently large $n$  the $\b$-iteration solves  Problem \ref{prob2}  \ (see Remark \ref{rem_5_1} below).

 \end{theorem}

In light of the fact that we only obtain approximate controls, we cannot expect to achieve exact asymptotic tracking. This leads us to the more reasonable problem referred to as Practical Asymptotic Regulation.
\begin{problem}[\bf  Practical  Regulation Problem]\label{prob2}
{\rm  \mbox{}\\
Given a desired approximation error level  $\e>0$,   find a
time dependent control law $u \in C_b(\bbr^+,\sU) $ such that 
\begin{equation}
  \limsup _ {t \to \infty} \| e(t) \|_{\sY} \leq \e,  \label{er1p}
\end{equation}
 for a given reference signal $r \in C_b^N(\bbr^+,\sY)$ and disturbance $d \in C_b^N(\bbr^+,\sD)$.}
\end{problem}

\begin{remark}\label{rem_5_1}{\rm
\begin{enumerate}
\item
It will be clear later that  if the reference and disturbance signals are not infinitely smooth, e.g., 
if they are in $C_b^N$ but not $C_b^{N+1}$,  then at this point our $\b$-iteration method will break down. This means, in that case, there  is a lower limit for the achievable value of the constant $\e$ in Problem \ref{prob2}. In other words, in order to guarantee that we can choose an arbitrarily small $\e$, using our algorithm,  we need the reference and disturbance signals to be in  $C_b^\i([0,\i))$.  Note that in the classical situation (cf. \cite{linreg}) the reference and disturbance signals are assumed to be generated by a finite dimensional, neutrally stable exo-system, in which case 
   the are always in  $C_b^\i([0,\i))$. 
   
   \item In the special case in which the  tracking and disturbance are harmonic signals  
  $$r = m_r+A_r \sin(\a_r t), \ \ \ d = m_d+A_d \sin(\a_d t),\ \ \ \a_r>0, \ \ \a_d>0,$$ the  bound on the   error $E_n(t)$ 
given in  \rr{mainest}   becomes 
$$  \limsup_{t>0}  \|r(t) -C(\ol{z}_n)(t)\|_\sY =  \limsup_{t>0}  \|E_n(t)\|_\sY \leq  D^n \lp |A_r| \a_r^n + |A_d| \frac{D_d}D \a_d^n\rp.$$
If $\ol{\a}=\max\{\a_r, \a_d\}$ and $\ol{\a} D<1$ then the $\b$-iteration solves  Problem \ref{prob2}, for sufficiently large $n$. On the other hand, 
simple numerical simulations for a one dimensional heat equation show that if $\ol{\a} D>1$ then the $\b$-iteration does not converge.

\end{enumerate}
}

\end{remark}

\subsection{Nonlinear Systems}\label{subsec_nonlinear}

For the nonlinear case, the situation is much more involved but some useful estimates are obtained in  \cite{aulisa2018analysis} for the tracking problem without disturbances.   
The starting point for analysis of the errors in the nonlinear case begins with the formulas \rr{eq2_a} and \rr{eq2_b} which we list again here
$$e_0(t) = G_0(t) + (K\ast r')(t) + (K_d\ast d')(t) + (H\ast F_0)(t),
$$
$$e_j(t)    = G_j(t)     + (K\ast e_{j-1}')(t) + (H\ast F_j)(t).$$
Here the norms of the operators $G_0(t)$ and $G_j(t)$ are  bounded above by functions $g_0(t)$ and $g_j(t)$ in the class $\calg$ of functions that are sums of products of polynomials times decaying exponentials.  Such functions decay exponentially to $0$ as $t$ goes to infinity. We also note that derivatives of such functions also decay exponentially.

Notice that the expressions in  \rr{eq2_a} and \rr{eq2_b} only differ from  \rr{eq2_a_fe0} and \rr{eq2_b_fe0} by the nonlinear terms
$(H\ast F_0)(t)$ and $(H\ast F_j)(t)$ where $H$ is defined in \rr{eq29_b}, $F_j$ in \rr{Fj} so that  
\begin{align}
   (H\ast F_0)(t)  & =  \int_0^t H  (t-\tau)     f (\ol{z}^0(\tau))  \, \dd\tau,  \label{H0} \\
   (H\ast F_j)(t) &  =  \int_0^t H  (t-\tau)     F_j (\tau)  \, \dd\tau. \label{Hj} 
\end{align}

The main difficulty in the analysis in the nonlinear case arises from the need to take time derivatives of the expressions which involves taking higher and higher order  Frech\'et derivatives for the nonlinear terms. 
To analyze the errors given in  \rr{eq2_a} and \rr{eq2_b}  we will need to estimate not only the norm of the nonlinear term $F_n$ but also derivatives of this nonlinear operator with respect to $t$. 

Here we assume that $f\, :\,  \sZ \ra \sZ $ is a nonlinear Lipschitz continuous function satisfying, for some $0<s_1<1$, and  for all $\vp_1$, $\vp_2 \in  \sZ$
  \begin{equation}
\label{eq_P}
\| f(\vp_1) - f(\vp_2)\|_{\sZ} \leq \e \|\vp_1-\vp_2\|_{s_1}.
\end{equation}
Therefore we have
\begin{equation}
\label{nFn}
\|F_n(t)\|_{\sZ} = \|f(\ol{z}_n )-f(\ol{z}_{n-1} )\|_{s_1}   \leq \e \|\ol{z}_n-\ol{z}_{n-1}\|_{s_1} = \e \|\ol{z}^n\|_{s_1} .
\end{equation}

According to Proposition 1, in  \cite{aulisa2018analysis}
\begin{proposition} \label{pr_bound_Fn}
\begin{equation}
\label{Fn}
\|F_n(t) \|_{\sZ^0} \leq \e  \|\ol{z}^n(t)\|_{s_1} .
\end{equation}
\end{proposition}
 Therefore, in order to estimate $F_n$ we need to estimate $\|\ol{z}^n(t)\|_{s_1}$.  

We will also need to estimate $\|F_n'(t)\|$ for which requires the Fr\'echet derivative of the Nemytskii operator for our nonlinear term. For the complete  details see  in \cite{aulisa2018analysis}.

We now turn to the estimate  for $e_0$. This corresponds to $E_1$ in \cite{aulisa2018analysis} which is studied in Section 5.1.  The main change that is needed  is that equation 5.6 in \cite{aulisa2018analysis} is replaced by 
\begin{equation}
\label{eq40}
e_0(t) = K\ast r'(t) + K_d \ast d'(t)  +H\ast F_0(t).  
\end{equation}

\begin{definition}\label{def1}{\rm
 For   $I \subseteq \bbr^+=\{t\, :\, 0\leq t<\i\}$   a fixed interval, the space $C_{b}^N(I)$ consisting of bounded and  $N$-times continuously differentiable  functions with derivatives bounded on $I$ is a Banach space  with the norm 
 $$\|\vp\|_{I,N} = \max_{0\leq j \leq N} \left(\sup_{t\in I}  |\vp^{(j)}(t)|\right).$$
 
 If $I=[0,\i)$ we write $\|\vp\|_{I,N} = \|\vp\|_{N}$. 
 

 }
 \end{definition}

Due to the presence of the disturbance this requires a slight modification of the analysis given in the proof of Proposition 4 in \cite{aulisa2018analysis} for the first step in the present analysis. Namely, with the obvious modifications to the proof of Proposition 4 we arrive at
\begin{proposition}\label{prop4}
\begin{equation}
\|e_0 \|_{ (0,t)}  \leq ( D     + \cald )  \|r\|_{1, (0,t)}+ ( D_d     + \cald_d )  \|d\|_{1, (0,t)}  ,\\ 
\end{equation}
and
\begin{equation}
\label{eq_e1_est}
\limsup \|e_{0}\|   \leq (D+\cald) \limsup \|r\|_1 +  ( D_d     + \cald_d )\limsup \|d\|_1  .  
\end{equation}
where $D$ and $D_d$ are defined in \rr{Dform} and 
\begin{align}
D_H =& \int_0^\i \| Ce^{A_\b t} I_0\| \, \dd t \ \   D_B =    \int_0^t \| e^{A_\b t} B  \| \, \dd t  \label{sK} \\
D_{A_\b} = & \int_0^t \| e^{A_\b t} I_0  \| \, \dd t \ \  D_{\Bd} =   \int_0^t \| e^{A_\b t} I_0 \Bd  \| \, \dd t \label{sP}\\
\cald = & D_H \lp \frac{D_B \e }{1- \e D_{A_\b}}\rp , \ \  
\cald_d =  D_H \lp \frac{D_{\Bd} \e }{1- \e D_{A_\b}}\rp,
\end{align}
where $\e$ is given in \rr{nFn}.
\end{proposition}
And,  for $n\geq 1$ we have 
\begin{proposition}
 \begin{align}
\label{en_2}
i) &\qquad\|e_n \|_{ (0,t)}  \leq ( D     + \cald )  \|e_{n-1}\|_{1, (0,t)}  ,\\
\mbox{ and }\qquad \nn &\\
\label{e1est}
ii) &\qquad \limsup \|e_{n}\|   \leq (D+\cald) \limsup \|e_{n-1}\|_1  .
\end{align}

\end{proposition}


 \begin{remark}\label{sinusodal}
   \begin{enumerate} 
  \item  Under our assumptions on $C$, $B$ and $B_d$ the values of $D$ and $D_d$ in \rr{Dform} are finite. In particular, for values of $s$ including $s_c$, $s_b$ and $s_d$  in Assumption \ref{ass1}, the scale of spaces $\sZ^{-s}$ generated by $A$ and the spaces $\sZ_\b^{-s}$ generated by $A_\b$ are equal  and  the norms in these spaces are equivalent (see, e.g.,   \cite{mikkola}, Lemma 9.4.3, page 443 or \cite{staffans}, Theorem 3.10.11, page 174).  Therefore  we can conclude that
 $C\in \sL(\sZ^{s_c},\sY)$, $B\in \sL(\sU,\sZ^{-s_b})$ and $D\in \sL(\sD,\sZ^{-s_d})$ if and only if $C\in \sL(\sZ_\b^{s_c},\sY)$, $B\in \sL(\sU,\sZ_\b^{-s_b})$ and $D\in \sL(\sD,\sZ_\b^{-s_d})$. In the following we denote the norm in $\sZ_\b^{s}  $ by $\|\cdot \|_{\b,s}$.

  \item In our estimates for the norms of $K(t)$  we will use  bounds in the scale of spaces $\sZ_\b^{-s}$. 
\begin{equation} \label{bddL1K}
  \|K(t) \|_\sL  \leq  \frac{ \|C\|_{\b, -s_c}    \|B\|_{\b, -s_b} \|\sab{-1+\d}\|  M_\b  \ e^{-\o_\b t}}{\b } \defeq \wt{Q} \,  e^{-\o_\b t},    
  \end{equation}
  so that 
\begin{equation}\label{eqD}
D  = \int_0^\i \|K(t)\|_{\call} \, \dd t \leq  \wt{Q} 
  \int_0^\i e^{-\o_\b t} \dd t  
   =   \frac{\wt{Q}  }{  \o_\b } \defeq D_H <\i.
  \end{equation}

\item For $K_d(t) = -CA^{-1}_{\beta} e^{A_\beta t} I_0 \Bd $  we  need   to estimate  $I_0\Bd$ where 
 $I_0 = (I+ BCA^{-1})$. Set $s_m = \max\{s_b, s_d\} $ and
\begin{align*}
 \| \sab{-s_m} \cald \| \leq&   \|\sab{-s_m}\Bd\|  + \| \sab{-s_m}B\|   \\
 &  \times  \| C\sa{-s_c}\| \| \sa{-1+s_c} \Bd\|  \defeq \wt{R}_{\cald} <\i.
 \end{align*}
 Notice that in the last  term we have $-1+s_c <-s_d$ due to our assumption $s_c+s_d <1$. Set $p = s_c+s_m<1$ and we have
\begin{align}
\|K_d(t)\|_{\call_d} &= \|  -C \sab{-1} e^{A_\beta t} ( I+ BCA^{-1})\Bd \| \nn  \\
  & \leq \|C\sab{-s_c}\| \ \| \sab{-1+p } \|\ \| e^{A_\beta t}  \|  \ \wt{R}_{\cald} \nn  \\
  & \leq   \|C\|_{\b,-s_c} \wt{R}_{\cald} \| \sab{-1+p} \|\   M_\b e^{-\o_\b t} = \wt{Q}_d e^{-\o_\b t}, \label{bddBd}
\end{align}
where we have defined
 \begin{equation}
\label{wtQd}
\wt{Q}_d =    \|C\|_{\b,-s_c} \wt{R}_{\cald}   M_\b \| \sab{-1+p} \|\   . 
\end{equation}
Then we have a  bound on the $L^1$ norm of $K_d(t)$  given by 
 \begin{equation}
\label{bddL1Kd}
 D_d =  \int_0^\i  \|K_d(t)\|_{\call_d }  \dd t \leq 
 \wt{Q}_d  \int_0^\i e^{-\o_\b t} \dd t  =
 \frac{\wt{Q}_d }{\o_\b} = D_{H_d} <\i . 
\end{equation} 
  
 \item For $0<s<1$, the definition of $\sab{-s}$ (see, for example, Definition 1.4.1 in \cite{henry})  gives
 \begin{equation}
\label{bddmAbs}
 \| \sab{-1+s} \| \leq \frac{M_\b}{\o_\b^{1-s}}. 
\end{equation}

 \end{enumerate} 

  \end{remark}

In the proof of  Proposition 4 in \cite{aulisa2018analysis} an  estimate of $\|\ol{z}^n(t)\|_{s_1}$ is required. We present the result of that estimate from \cite{aulisa2018analysis} here
\begin{proposition} \label {pr_limsup_zn} 
 \begin{align}
\label{e_n_t}
 i) & \qquad \|\ol{z}^n \|_{s_1,(0,t)} \leq \frac{D_B}{1-\e_1  D_{A_\b}} \|e_{n-1}\|_{(0,t)},
\\
\mbox{ and }\qquad \nn&
\\
\label{limsup_zn*}
ii) & \qquad \limsup   \|\ol{z}^n \|_{s_1} \leq \frac{D_B}{1-\e_1  D_{A_\b}}   \  \limsup  \| e_{n-1} \| .
\end{align}
\end{proposition}

Similar estimates for $e_2(t)$ are obtained in  Section 5.2 in \cite{aulisa2018analysis}. We have from \rr{eq2_b} 
$$e_2(t) = K\ast e_1'(t) +H\ast F_2(t). $$
As pointed out in  \cite{aulisa2018analysis}, the main technical difficulty here comes from the fact that we need to differentiate  $e_0(t) $ which gives
\begin{align*}
e_0'(t) =&  (K\ast r')' (t) +(K_d\ast d')' + (H\ast F_1)'  \\
=&  (K\ast r'')(t) + K(t)r'(0)+  (K_d\ast d'')(t) + K_d(t)d'(0)+ (H\ast F_0')(t) + H(t) F_0(0).
\end{align*}
Since $F_n(0)=0$ this becomes
 \begin{equation}
\label{e1p}
e_0'(t) =  (K\ast r'')(t) + K(t)r'(0)+(K_d\ast d'')(t) + K_d(t)d'(0) + (H\ast F_0')(t) .
\end{equation} 
Notice that
 both $K(t)r'(0) $ and $K_d(t)d'(0)$ are terms whose norms give functions in the class $\calg$ and therefore decay to zero exponentially in time. 
So clearly the main new ingredient is $ F_1'$.

Using a modification  of the arguments given in Section 5.2 in  \cite{aulisa2018analysis} we obtain the analog of Theorem 4 in \cite{aulisa2018analysis}
\begin{theorem}
\begin{equation}
\label{e2_est}
\limsup \|e_1\| \leq  (D+\cald)^2 \left[ \limsup \|r\|_2 + \frac{(D_d+\cald_d)}{(D+\cald)} \limsup \|d\|_2\right]. 
\end{equation}
\end{theorem}

For $n>1$ the calculations become much more involved and we were not able to repeat a simple estimate like the one in \rr{e2_est}.  For the third iteration, i.e., for $n=2$,  we do obtain 
$$\limsup \|e_2\| \leq  (D+\cald)^3 \left[ \limsup \|r\|_3 + \frac{(D_d+\cald_d)}{(D+\cald)}  \limsup \|d\|_3\right] +V_2,$$
 where $V_2$ is a complicated expression which is usually very small.

 \end{document}

%% file: dsgheader20.tex
\usepackage{amsmath,latexsym,amssymb}

\def\tra{{\intercal}}

\def\nn{\nonumber}

\def\dd{\mathrm d}
\def\exp{\mathrm e}

 \def\sb{s_b}
  \def\sd{s_d}
   
   \def\sc{s_c}
   
    \def\Bin{B_{\text{in}}}
  \def\Bd{B_{\text{d}}}

\newcommand{\szmb}{\sZ^{-s_b}}

\newcommand{\szmd}{\sZ^{-s_d}}

\newcommand{\sa}[1]{(-A)^{#1}}
\newcommand{\sab}[1]{(-A_\b)^{#1} }

\newcommand{\eqdef}{=\mathrel{\mathop:}}
\newcommand{\defeq}{ \mathrel{\mathop:}=}

\newcommand{\deriv}[3][]{%
  \ensuremath{\mathchoice{\frac{d^{#1}{#2}}{d{#3}^{#1}}}%
  {{d^{#1}{#2}/d{#3}^{#1}}}
  {{d^{#1}{#2}/d{#3}^{#1}}}%
  {{d^{#1}{#2}/d{#3}^{#1}}}}}
\newcommand{\pderiv}[3][]{%
  \ensuremath{\mathchoice{\frac{\partial^{#1}{#2}}{\partial{#3}^{#1}}}%
  {{\partial^{#1}{#2}/\partial{#3}^{#1}}}%
  {{\partial^{#1}{#2}/\partial{#3}^{#1}}}%
  {{\partial^{#1}{#2}/\partial{#3}^{#1}}}}}

\newcommand{\rr}[1]{\eqref{#1}}

\def\nn{\nonumber}

\newcommand{\lp}{\left(}
\newcommand{\rp}{\right)}

\usepackage{graphicx,euscript}

\newcommand{\script}[1]{\EuScript{#1}}
\newcommand{\sA}{{\script A}}
\newcommand{\sB}{{\script B}}

\newcommand{\sD}{{\script D}}

\newcommand{\sF}{{\script F}}

\newcommand{\sL}{{\script L}}
\newcommand{\sM}{{\script M}}

\newcommand{\sP}{{\script P}}

\newcommand{\sT}{{\script T}}
\newcommand{\sU}{{\script U}}

\newcommand{\sW}{{\script W}}
\newcommand{\sX}{{\script X}}
\newcommand{\sY}{{\script Y}}
\newcommand{\sZ}{{\script Z}}

\usepackage{graphicx}

\newcommand{\beqn}{\begin{eqnarray*}}
\newcommand{\eeqn}{\end{eqnarray*}}

\def\calg{\mathcal{G}}

\def\call{\mathcal{L}}

\def\cald{\mathcal{D}}
\def\calg{\mathcal{G}}

\DeclareMathOperator{\ran}{Ran}

 \def\z{\zeta}
 \def\o{\Omega}
\def\g{\gamma}
\def\G{\Gamma}
\def\e{\epsilon}
\def\a{\alpha}
\def\b{\beta} 
\def\d{\delta}
\def\r{\rho}

\def\s{\sigma}
\def\i{\infty}
\def\l{\lambda}

\def\o{\omega}
\def\ra{\rightarrow}
\def\vp{\varphi}

\def\wh{\widehat}
\def\wt{\widetilde}

\def\i{\infty}
\def\l{\lambda}

\def\o{\omega}

\def\ra{\rightarrow}

\def\wh{\widehat}
\def\wt{\widetilde}

\def\bbc{\mathbb{C}}
\def\bbz{\mathbb{Z}}

\def\bbr{\mathbb{R}}

\newcommand{\ol}{\overline}

\def\re{\text{Re }}